\documentclass[proceedings]{aofa}
\usepackage{amsmath}

\usepackage{amsthm}

\title{Variance of additive functions defined on random assemblies}

\author{ Eugenijus Manstavi\v cius and Vytautas Stepas}

\address{Department of Mathematics and Informatics, Vilnius University, Naugarduko str. 24, LT-03225 Vilnius,  Lithuania}

\keywords {Labeled decomposable structure, additive function, moments, Tur\'an-Kubilius inequality.}

\begin{document}

\maketitle

\begin{abstract}
\paragraph{Abstract.}
An inequality for the variance of an additive function defined on random decomposable structures, called assemblies, is established. The result generalizes estimates obtained earlier in the cases of  permutations and mappings of a finite set into itself. It is analogous to the Tur\'an-Kubilius inequality for additive number-theoretic  functions.

\end{abstract}

\newtheorem{thm}{Theorem}
\newtheorem{lem}{Lemma}
\newtheorem*{klem}{Key Lemma}
\newtheorem{cor}{Corollary}
\newtheorem*{cor*}{Corollary}
\newtheorem{conj}{Conjecture}
\newtheorem*{def*}{Definition}
\newtheorem*{cond*}{Condition 1}

\newtheorem{prop}{Proposition}
\newtheorem*{prop*}{Proposition}
\newtheorem*{CLT}{CLT (\cite{EM-LMJ96})}
\newtheorem*{Comp}{Compactness Thm (\cite{EM-RJ08})}

\newtheorem*{WLLN}{WLLN (\cite{EM-RJ08})}

\def\E{\mathbf{E}}

\def\R{\mathbf{R}}
\def\Z{\mathbf{Z}}
\def\N{\mathbf{N}}
\def\C{\mathbf{C}}
\def\D{\mathbf{D}}
\def\G{\mathbf{G}}
\def\S{{\mathbf{S}_n}}
\def\E{\mathbf{E}}
\def\V{\mathbf{V}}
\def\cG{\mathcal{G}}

\def\cV{\mathcal{V}}
\def\cA{\mathcal A}
\def\cU{\mathcal U}
\def\cF{\mathcal F}

\def\Ra{\Rightarrow}

\def\k{\kappa}

\def\e{\varepsilon}
\def\ro{{\rm o}}

\def\rO{{\rm O}}
\def\re{{\rm e}}
\def\rd{{\rm d}}

\def\oJ{\overline{J}}
\def\n{$n\to\infty$}

\def\cF{\mathcal F}
\def\s{\smallskip}

    We deal with additive functions defined on combinatorial structures such as permutations, mappings of a finite set into itself, 2-regular graphs etc. If a structure is taken at random, such functions are sums of dependent random variables; sometimes, they are called \textit{separable statistics}. Their value distribution is a complex problem.  One of the useful tools in analysing it are estimates of the variance. This is our main objective. On the other hand, our interest has been highly stimulated by the Tur\'an-Kubilius inequality
in probabilistic number theory or by analogous inequalities in the theory of additive arithmetical semigroups.

   An assembly is a construction defined on a set  by its partition and some structure introduced in all subsets, afterwards called components of the assembly. Assume that given a subset of size $j$ we can  introduce  $g_j<\infty$ structures, then the number of assemblies spanned over an $n$ set (assemblies of the order $n$)  equals
    \[
       G(n)=n!\sum_{\ell(\bar s)=n}\prod_{j=1}^n\Big(\frac{g_j}{j!}\Big)^{s_j}\frac{1}{s_j!}=:n!Q(n).
    \]
    Here $n\in\N$, $\ell(\bar s):=1s_1+\cdots+ns_n$ if $\bar s=(s_1,\dots,s_n)\in \N_0^n$ and the summation is over such vectors satisfying $\ell(\bar s)=n$. We will denote the class of assemblies by $\cG$ and the set of assemblies of the order $n$ by ${\cG}_n\subset \cG$.
    
     Let $\lambda_j:=g_j/j!$. In the past decades much attention was paid to the \textit{logarithmic} class defined by the asymptotic condition
$\rho^j j\lambda_j\sim\theta$ for some positive constants $\theta$, $\Theta$ and $\rho$ as $j\to\infty$ (see \cite{ABT}). Extensions were initiated in  the first author's paper \cite{EM-CPC02}, where
a condition
\begin{equation}
         0<\theta\leq  \rho^j j\lambda_j\leq \Theta,\quad j\geq 1,
\label{cC}
\end{equation}
 was used. The lower bound excluded, for example, the class of 2-regular graphs, however. Basing upon the experience, in the present paper we confine ourselves to a class of assemblies characterized by
 some positive constants $\rho$, $\Theta$, $\theta$, $\theta'$, and  $n_0\geq 1$.

 \begin{def*} We say that a class of assemblies is weakly logarithmic if the following conditions are satisfied:
\begin{equation}
\label{111}
\rho^jj\lambda_j\leq\Theta, \quad j\geq 1;
\end{equation}
\begin{equation}
\label{121}
\sum_{j\leq n}\rho^jj\lambda_j\geq \theta n,\quad n\geq n_0;
\end{equation}
\begin{equation}
\label{131}
n Q(n)\rho^n\geq \theta' \exp\bigg\{\sum_{j\leq n}\lambda_j \rho^j\bigg\}, \quad  n\geq 1.
\end{equation}
\end{def*}

Let $k_j(\sigma)\geq 0$ be the number of components  of size  $j$ in $\sigma\in \cG_n$ and $1\leq j\leq n$. An \textit{additive function} $h:\cG_n\to\R$ is  defined by a
real two-dimensional array $\{h_j(k)\}$, where $j,k\in\N$,  $jk\leq n$,
 and $h_j(0)=0$ for all $j\leq n$,  by setting
 \begin{equation}
                 h(\sigma)=\sum_{j=1}^n h_j\big(k_j(\sigma)\big).
 \label{h}
 \end{equation}
  Apart from the most popular example of the number-of-components function $w(\sigma)=k_1(\sigma)+\cdots+k_n(\sigma)$, they appear in many algebraic and combinatorial problems. Particular additive functions appear in physical models as a part of Hamiltonians in the Bose gas  theory.
  
Let $\E_nh$ and ${\V}_nh$ denote the expectation and the variance of $h=h(\sigma)$ with respect to uniform probability measure. The problem is to estimate
\[
  {\V}_n h=\frac{1}{G(n)}\sum_{\sigma\in \cG_n}\left(h(\sigma)-\E_n h\right)^2=\E_n h^2-(\E_n h)^2
  \]
  in terms of the values $h_j(k)$ where $jk\leq n$ and parameters characterizing the class of assemblies.
 
In the sequel, let $Q^{\{j\}}(n)$ be defined by 
\[
Q^{\{j\}}(n)=\sum_{\ell(\bar s)=n\atop s_j=0} \prod_{i\leq n}{\lambda_i^{s_i}\over s_i!},\quad j\leq n
\]
and $\ll$ be an analog of the symbol $O(\cdot)$.

Now, the results.

\begin{thm} \label{thm1} Assume that $\cG$ is weakly logarithmic and  $h:\cG_n\to\R$ is an arbitrary additive function. Then
\begin{equation}
             {\mathbf V}_n h\ll \sum_{jk\leq n}{\lambda_j^kh_j(k)^2\over k!} {Q^{\{j\}}(n-jk)\over Q(n)}
             \label{result}
             \end{equation}
for $n\geq 1$.
\end{thm}

Inequality (\ref{result})  sharpens a bit Theorem 3 in \cite{EMVS-AofA14} proved for an arbitrary  additive function defined on weighted permutations under  condition (\ref{cC}).

A \textit{completely additive function} $h$ is defined by the array $h_j(k)=a_jk$, where $a_j\in\R$ and $jk\leq n$.
For such functions, inequality (\ref{result}) takes a simpler form.

\begin{thm} \label{thm2} Assume that $\cG$ is weakly logarithmic and  $h:\cG_n\to\R$ is a completely additive function  defined via $h_j(s)=a_j s$ where $js\leq n$. Then
\begin{equation}
             {\mathbf V}_n h\ll \sum_{j\leq n}\lambda_j a_j^2 {Q(n-j)\over Q(n)}.
             \label{result1}
             \end{equation}
\end{thm}

Inequality (\ref{result1}) for weighted permutations satisfying  condition (\ref{cC}) has a longer history (see \cite{EMVS-AofA14})).

\end{document}